\documentclass[11pt]{amsart}
\usepackage{
  amsmath,
  amsfonts,
  amssymb,
  graphicx, 
  times,
  color,
  hyphenat,
  mathptmx,
  pinlabel
}
\input xy
\xyoption{all}
\newcommand{\rh}{\overline{H}}
\newcommand{\trh}{\overline{\overline{H}}}
\newcommand{\rhn}{\overline{H}_N}
\newcommand{\trhn}{\overline{\overline{H}}_N}

\newcommand{\cH}{\mathcal{H}}
\newcommand{\cP}{\mathcal{P}}

\newcommand{\n}{\noindent}

\newcommand{\figins}[3] 
{\raisebox{#1pt}{\includegraphics[height=#2 in]{figs/#3}}}

\newtheorem{thm}{Theorem}[section]
\newtheorem{lem}[thm]{Lemma}

\theoremstyle{definition}

\newcommand{\bQ}{\mathbb{Q}}

\newcommand{\cT}{\mathcal{T}}

\newcommand{\ra}{\rightarrow}
\newcommand{\xra}[1]{\xrightarrow{#1}}
\newcommand{\lra}{\longrightarrow}

\catcode`\@=11
\long\def\@makecaption#1#2{%
    \vskip 10pt
    \setbox\@tempboxa\hbox{%
\small{#1: }\ignorespaces #2}%
    \ifdim \wd\@tempboxa >\captionwidth {%
        \rightskip=\@captionmargin\leftskip=\@captionmargin
        \unhbox\@tempboxa\par}%
      \else
        \hbox to\hsize{\hfil\box\@tempboxa\hfil}%
    \fi}
\newdimen\@captionmargin\@captionmargin=2\parindent
\newdimen\captionwidth\captionwidth=\hsize
\catcode`\@=12
\title{The reduced HOMFLY-PT homology for the Conway and the Kinoshita-Terasaka knots}
\author{Marco Mackaay}
\address{Departamento de Matem\'{a}tica\\ Universidade do Algarve\\ 
Campus de Gambelas\\ 8005-139 Faro\\ Portugal and CAMGSD\\Instituto Superior T\'{e}cnico\\ Avenida Rovisco Pais\\ 
1049-001 Lisboa\\ Portugal}
\email{mmackaay@ualg.pt}
\author{Pedro Vaz}
\address{Departamento de Matem\'{a}tica\\ Universidade do Algarve\\ 
Campus de Gambelas\\ 8005-139 Faro\\ Portugal and  
CAMGSD\\Instituto Superior T\'{e}cnico\\ Avenida Rovisco Pais\\ 
1049-001 Lisboa\\ Portugal}
\email{pfortevaz@ualg.pt}
\begin{document}
%
%
\newdimen\captionwidth\captionwidth=\hsize
%
%
\begin{abstract} 
In this paper we compute the reduced HOMFLY-PT homologies of the 
Conway and the Kinoshita-Terasaka knots and show that they are isomorphic. 
\end{abstract}
\maketitle
%
%
%
\section{Introduction}\label{sec:intro}

In this paper we use Rasmussen's results in~\cite{rasmussen-diff} to compute the 
reduced HOMFLY-PT homologies, defined by Khovanov and Rozansky~\cite{KR2}, of the 
Conway and the Kinoshita-Terasaka knots. It turns out that these homologies are isomorphic. We also show that our calculations imply that the 
Khovanov-Rozansky $sl(N)$-homologies of these two knots 
are isomorphic for all $N\geq 2$.  
This result surprised us 
because the Floer knot homologies of these knots are non-isomorphic~\cite{OS}. 
Since people have conjectured that for each knot there should exist a spectral 
sequence converging to the Floer knot homology with $E_2$-page isomorphic to 
the HOMFLY-PT homology, our result shows that the differentials of 
the conjectured spectral sequences for the Kinoshito-Terasaka and the Conway 
knot should be different. 

We did 
our calculations in the summer of 2006 and simply put them in a drawer. Since then 
several people, who knew about the result, asked us to write it up, which 
is why we finally decided to write this small note.

We claim no original theoretical insights. We have simply used Rasmussen's 
results. Since these calculations are hard and we had to use 
all sorts of tricks, this note might help other people to understand Rasmussen's results and to do calculations themselves and it also gives the reduced HOMFLY-PT homologies of the aforementioned knots, 
which to our knowledge were not available before.

\section{Rasmussen's toolkit}
\label{sec:rasmussen}

In this section we explain Rasmussen's results~\cite{rasmussen-diff} which 
enable us to do the computations. Note that we take the usual triple gradings 
$(a,q,t)$ from the HOMFLY-PT homology, using the conventions in \cite{KR2}, rather 
than Rasmussen's gradings in~\cite{rasmussen-diff}. We thank Rasmussen for explaining the 
conversion rules between the two gradings. 

Given a knot $K$, denote its reduced HOMFLY-PT homology by $\rh(K)$ and 
its reduced $sl(N)$ homology by $\rhn(K)$.
\begin{thm}\emph{(Thm. 2 in~\cite{rasmussen-diff})} For each $N>0$ there is a spectral sequence 
$\bigl(E_k(N), d_k(N)\bigr)$ which starts at $\rh(K)$ and converges to 
$\rhn(K)$.  
\label{thm:ssN}
\end{thm}
\n The $(a,q,t)$-degree of the $d_k(N)$ is $(-2k,2Nk,1)$. Note that for $N$ big enough the differentials $d_k$ are zero.

\begin{thm}\emph{(Thm. 3 in~\cite{rasmussen-diff})} There is a spectral sequence $\bigl(E_k(-1),d_k(-1)\bigr)$ starting 
at $\rh(K)$ and converging to $\bQ$. 
\label{thm:ss1}
\end{thm}
\n The $(a,q,t)$-degree of $d_k(-1)$ is $(2-2k,2-2k,2k-1)$.
Another result that will be useful for us is the following.
\begin{lem}\emph{(Lemma 6.2 in~\cite{rasmussen-diff} )}
\label{lem:acom}
For any $k$ the differentials $d_k(-1)$ and $d_1(N)$ anticommute.
\end{lem}

For a two-component link $L$ denote by $\rhn(L,i)$ the $sl(N)$ homology of 
$L$ reduced w.r.t. to the link component $i$. For $j\ne i$, let 
$$\rhn(L,i){\xra{\ \ X_j\ \ }}\rhn(L,i)$$ 
be the map induced by multiplication by $X_j$ (i.e. on the $j-th$ component of $L$). 
Note that this map has $(q,t)$-bidegree $(2,0)$. Let $\trhn(L)$ 
be the totally reduced $sl(N)$ homology of $L$.

\begin{lem} There is a long exact sequence 
\begin{equation}
\label{eq:totred}
\dotsi\lra\rhn(L,i)\ {\xra{\ \ {X_j}\ \ }}\
\rhn(L,i)\ {\xra{\ \ {(-1,1/2)}\ \ }}\ \trhn(L)\
{\xra{\ \ {(-1,1/2)\ \ }}}\ \rhn(L,i)\lra\dotsi
\end{equation}
\end{lem}  
We have not given all the maps explicitly, since we do not need them, but we have 
given their bidegrees. 

Let $K_+$ be knot with a given positive crossing. Let $K_-$ be the same knot except 
for that particular crossing which is now negative and let $K_0$ be the two-component 
link obtained from $K$ by the oriented resolution of the same crossing. 
\begin{lem}\emph{(Lemma 7.6 in~\cite{rasmussen-diff} )}
 There is a long exact sequence 
\begin{equation}
\label{eq:les}
\dotsi\lra\rhn(K_-)\
{\xra{\ \ {(N,-1/2)}\ \ }}\
\trhn(K_0)\ {\xra{\ \ {(N,-1/2)}\ \ }}\ \rhn(K_+)\
{\xra{\ \ {(-2N,2)\ \ }}}\ \rhn(K_-)\lra
\dotsi
\end{equation}
\end{lem}  
  
Let us do a simple example to illustrate the exact sequences above. By abuse of notation we always identify the homology with the Poincar\'e polynomial. Thus, multiplying the homology by a polynomial means multiplying the Poincar\'e polynomial by that polynomial.


Recall (see \cite{rasmussen-diff}) that the positive Hopf link, $\cH^+$, has reduced 
HOMFLY-PT 
homology equal to 
$$aq^{-1}+q(q^Nt^{-1})^2(q^{-N+2}+q^{-N+4} + \dotsb + q^{N-2}).$$
Note that it does not matter which component we choose for reduction in this case. 
Multiplication by $X$ maps the generator of degree $q^{s}$ to the generator of 
degree $q^{s+2}$ in $q^{-N+2}+q^{-N+4} + \dotsb + q^{N-2}$. Therefore the kernel of 
$$\rhn(\cH^+,i){\xra{\ \ X_j\ \ }}\rhn(\cH^+,i)$$ 
is given by the generators $aq^{-1}+q^{N-1}(q^Nt^{-1})^2$ and the 
cokernel by the generators $aq^{-1}+q^{-N+3}(q^Nt^{-1})^2$. Using the exact sequence 
\eqref{eq:totred} we see that 
$$
\trh(\cH^+)=a^3t^{-5/2}+aq^{-2}t^{1/2}+at^{-1/2}+aq^2t^{-3/2}.
$$   
Let $K_+$ be the positive trefoil $\cT^+$, then (see~\cite{rasmussen-diff})
$$\rh(K_+)=a^2q^{-2}+a^2q^2t^{-2}+a^4t^{-3}.
$$
Note that $K_0=\cH^+$ and $K_-$ is the unknot. One now easily checks that the long 
exact sequence \eqref{eq:les} holds.

There is also a useful variant of the exact sequence~\eqref{eq:totred}. Let $L_-$ and $L_+$ be the two-component link diagrams which differ by the sign of one crossing between different components and $K$ the knot diagram obtained by resolving that crossing respecting the orientations.
\begin{lem} 
There is a long exact sequence 
\begin{equation}
\label{eq:ktotred}
\dotsi\lra\rhn(L_-)\ {\xra{\ \  (N,-1/2) \ \ }}\
\trhn(K)\ {\xra{\ \ {(N,-1/2)}\ \ }}\
\rhn(L_+)\ {\xra{\ \ {(-2N , 2)\ \ }}}\
\rhn(L_-)\lra\dotsi
\end{equation}
\end{lem}  

Notice that $X$ acts as zero on $\rh(K)$ and
\begin{equation}
\label{eq:ktred}
\trhn(K)=Cone\bigl(\rhn(K) \xra{\ \ 0 \ \ } \rhn(K)\bigr)
=\rhn(K)(q^{-1}t^{1/2}+ qt^{-1/2}).
\end{equation}

\section{Computations}\label{sec:comp}

In this section we compute the HOMFLY-PT homologies of the Conway and the Kinoshita-Terasaka knots. Resolving and changing a particular crossing of a diagram results in simpler diagrams where homology can be computed more easily. We then use the long exact sequence~\eqref{eq:les} recursively to obtain the homology of the initial diagram. The spectral sequences of Theorems~\ref{thm:ssN} and~\ref{thm:ss1} and the knowledge of the $sl(N)$ homology of the initial diagram help us to unambiguously identify isomorphisms and zero maps in the long exact sequences~\eqref{eq:les} and~\eqref{eq:ktotred}.

The HOMFLY-PT polynomial of the Conway and Kinoshita-Terasaka knots is
\begin{gather*}
\begin{split}
\cP(K_\text{Conway})= \cP(K_\text{KT})=
a^{-4}(q^{-4} - q^{-2} +2 - q^2 + q^4) + a^{-2}( -q^{-6} -2q^{-2} -2q^2 -q^6) 
\\
+ (q^{-6} +2q^{-2} +1 +2q^2 +q^6) + a^2( -q^{-4} +q^{-2} -2 +q^2 -q^4),
\end{split}
\end{gather*}
and their reduced Khovanov homologies are isomorphic and given by
\begin{gather*}
\label{eq:Kh}
\begin{split}
\overline{Kh}(K_\text{Conway}) &= Kh(K_\text{KT})=  q^8t^{-5} + 2q^6t^{-4} + 2q^4t^{-3}  + 3q^2t^{-2} + (3 + q^2)t^{-1}  + 3 + 2q^{-2}
\\  & \mspace{20mu}    
+ (2q^{-2} + 2q^{-4})t + (3q^{-4} + q^{-6})t^2  + 3q^{-6}t^3 + 2q^{-8}t^4 + 2q^{-10}t^5 + q^{-12}t^6 .
\end{split}
\end{gather*}
From Theorem~\ref{thm:ssN} it follows that $\dim{\rh(K_\text{Conway})}\geq 33$ and $\dim{\rh(K_\text{KT})}\geq 33$.

\subsection{The Kinoshita-Terasaka knot}\label{ssec:KT}

\begin{figure}[h!]
\centering
\figins{0}{1.8}{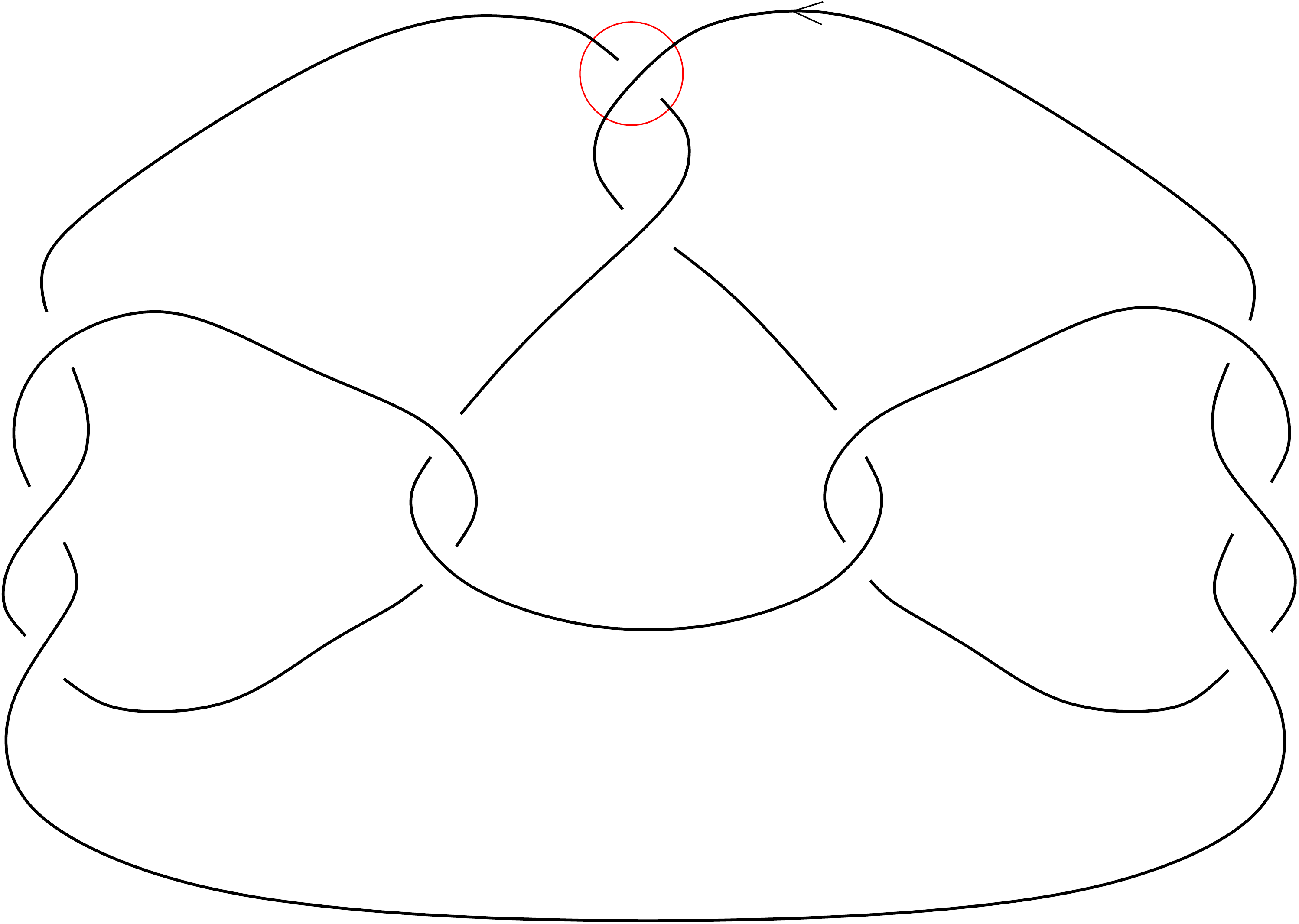}
\caption{The Kinoshita-Terasaka knot}
\label{fig:KT}
\end{figure}

\begin{figure}[h!]
\centering
\figins{0}{1.8}{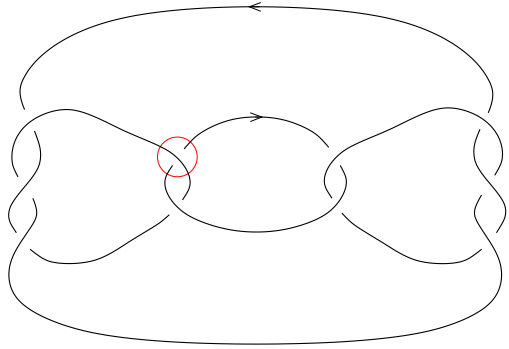}
\caption{Diagram $K_0$ obtained by taking the oriented resolution of the encircled crossing of the diagram of the KT knot of Figure~\ref{fig:KT}. The (negative) encircled crossing is used for the long exact sequence~\eqref{eq:ktotred}}
\label{fig:ktL0}
\end{figure}

A diagram of the KT knot is given in Figure~\ref{fig:KT}. 
Changing the encircled crossing we get the unknot, while resolving it results in the two-component link diagram $K_0$  depicted in Figure~\ref{fig:ktL0}.
Notice that $K_0$ is the pretzel link $P(3,-2,2,-3)$ which is amphicheiral ($K_0=K_0^!$). This diagram corresponds to the (non-alternating) link $L_{10n36}$ in Thistlethwaite's link table~\cite{katlas}.

\begin{figure}[h!]
\hair 2pt
\labellist
\pinlabel $M_0$ at  435 -50
\pinlabel $M_+$ at 1420 -50
\endlabellist
\centering
\figins{0}{1.75}{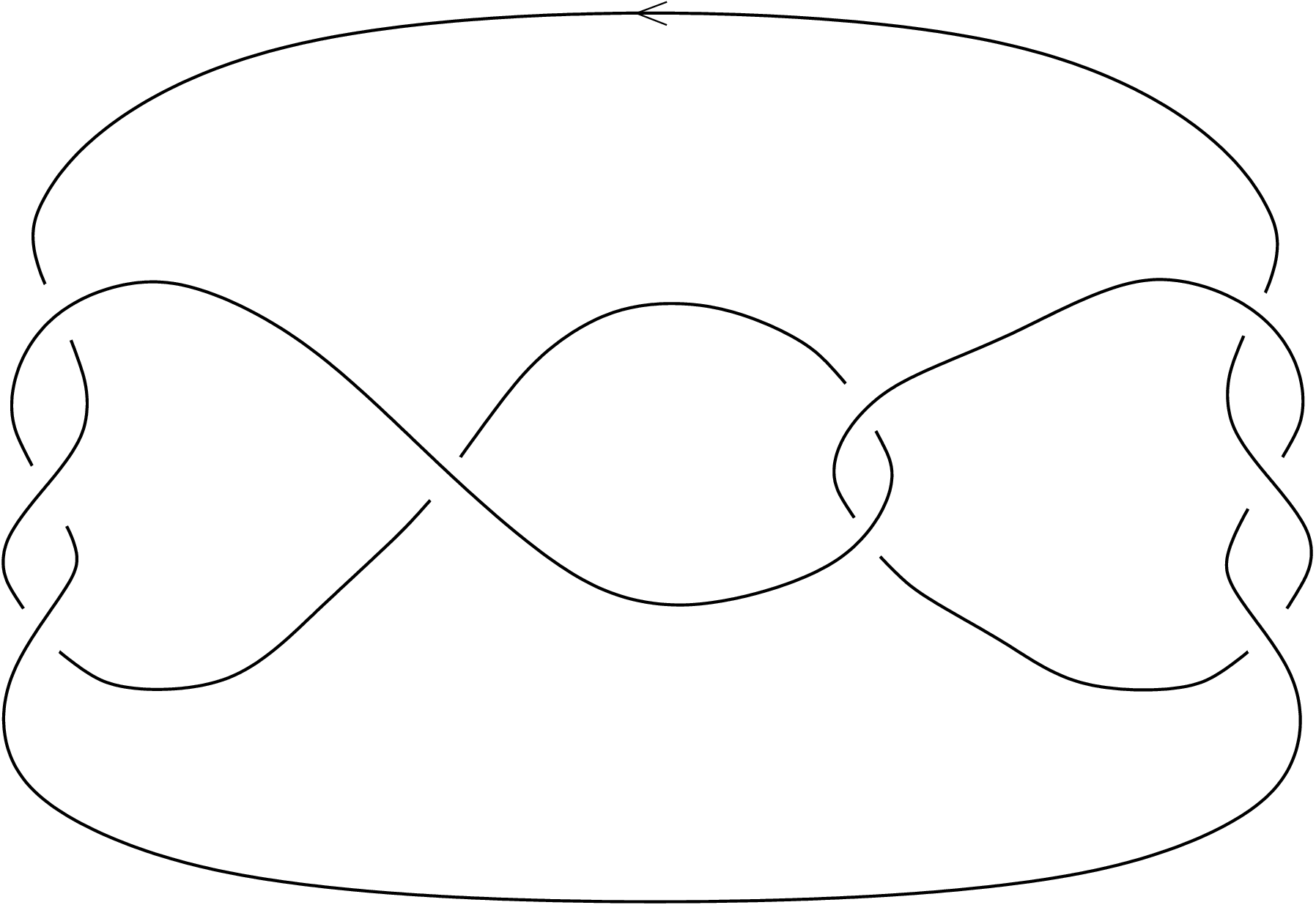}\qquad\ \figins{0}{1.75}{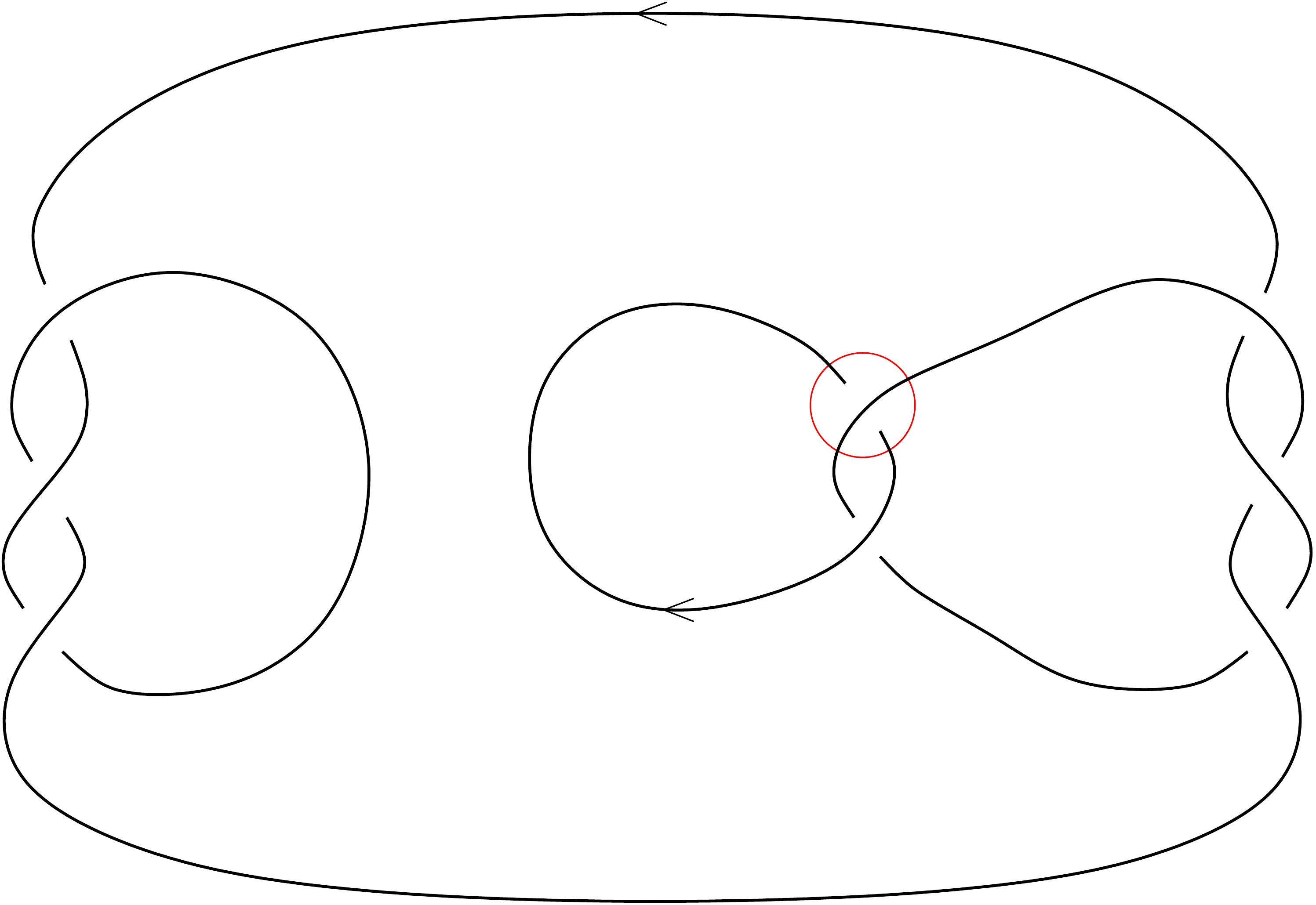}\\[2ex]
\caption{Diagrams $M_0$ and $M_+$ obtained by taking the oriented resolution of the encircled crossing of the diagram $K_0$ of Figure~\ref{fig:ktL0} }
\label{fig:ktM0}
\end{figure}

Taking the oriented resolution and changing the encircled crossing of $K_0$ we obtain the diagrams $M_0$ and $M_+$ respectively, depicted in Figure~\ref{fig:ktM0}.
The diagram $M_0$ corresponds to the pretzel knot $P(3,-1,2,-3)$ and is isotopic to the Pretzel knot $P(3,-2,-3)$, which in turn is the mirror of the knot $8_{20}$ in Rolfsen table~\cite{katlas}. We have 
\begin{gather}
\label{eq:HH820m}
\begin{split}
\trh(M_0) &= \trh(8_{20}^!) =
   a^{4}q^{3}t^{-11/2} 
+ (a^{2}q^{5} + a^{4}q)t^{-9/2} 
+ (a^{4}q^{-1} + a^{2}q^{3})t^{-7/2}
\\ & \mspace{20mu} 
+ (2a^{2}q +  a^{4}q^{-3})t^{-5/2} 
+  (q^{3}+2a^{2}q^{-1})t^{-3/2}
+ (2q+a^{2}q^{-3})t^{-1/2} 
\\ & \mspace{20mu} 
+ (2q^{-1} + a^{2}q^{-5})t^{1/2} 
+ q^{-3}t^{3/2}
.
\end{split}
\end{gather}

The diagram $M_+$ is the connected sum of the (positive) Hopf link with the connected sum of the positive trefoil $\cT^+$ and the negative trefoil $\cT^-$. From Lemma 7.8 of~\cite{rasmussen-diff} for connected sums it follows that

\begin{gather}
\label{eq:Hktm+}
\begin{split}
\rh(M_+) &= \rh(\cH^+)\otimes\rh(\cT^+)\otimes\rh(\cT^-)
\\ &=
   a^{3}qt^{-3} + aq^{3}t^{-2} + a^3q^{-3}t^{-1} + 3aq^{-1} + a^{-1}qt + aq^{-5}t^2 + a^{-1}q^{-3}t^3
\\ & \mspace{20mu}
+ [N-1]\bigl(
a^4q^3t^{-5} + a^2q^5t^{-4} + a^4q^{-1}t^{-3} + 3a^2qt^{-2} + q^3t^{-1}+a^2q^{-3}+q^{-1}t
\bigr)
.
\end{split}
\end{gather}

We can now determine $\rh(K_0)$ using the long exact sequence~\eqref{eq:ktotred} where the diagrams $K_0$, $M_0$ and $M_+$ correspond to $L_-$, $K$ and $L_+$ respectively. It reads
$$\dotsm\xra{\quad } \rh(K_0) \xra{\ at^{-1/2}}\trh(M_0)\xra{\ at^{-1/2}}\rh(M_+)\xra{\ a^{-2}t^2}\rh(K_0)\xra{\quad }\dotsm$$

From Equations~\eqref{eq:HH820m} and~\eqref{eq:Hktm+} and comparing degrees we have that
$t^{-11/2}a^4q^3$, $t^{-7/2}a^4q^{-1}$ and $t^{-1/2}q$
are in the kernel of the map $\trh(M_0)\xra{\ at^{-1/2}}\rh(M_+)$ and that 
$t^{-5}a^4q^2(q^{N-3}+\dotsm + q^{5-N})$, 
$t^{-4}a^2q^4(q^{N-3}+\dotsm + q^{3-N})$,
$t^{-3}a^4q^{-2}(q^{N-3}+\dots + q^{5-N})$,
$t^{-2}3a^2(q^{N-3}+\dotsm + q^{3-N})$, $t^{-2}a^2q^{N-1}$,
$t^{-1}q^2(q^{N-3}+\dotsm + q^{3-N})$,
$t^0a^2q^{-4}(q^{N-3}+\dotsm +q^{3-N})$, $t^0aq^{-1}$,
$tq^{-2}(q^{N-3}+\dotsm +q^{3-N})$, $ta^{-1}q$,
$t^2aq^{-5}$ and $t^3a^{-1}q^{-3}$
are in the cokernel.
Using this can we form a first list of guaranteed and possible generators of $\rh(K_0)$. Since $K_0^!=K_0$ the polynomial of $\rh(K_0)$ has to be invariant under the transformation $\psi(a,q,t)=(a^{-1},q^{-1},t^{-1})$. To have this symmetry we need to promote some possible generators to generators of $\rh(K_0)$ and discard possible generators not paired by $\psi$. Then we apply the exact sequence~\eqref{eq:ktotred} again to the newly promoted generators to obtain a new list which is in Table~\ref{tab:1st-KT}.
\begin{table}[h!]
\caption{Guaranteed and possible generators of $\rh(K_0)$}
\label{tab:1st-KT}
\newcommand\s[1]{\small{$#1$}}
\begin{center}
\begin{tabular}{|c|c|c|c|}
\hline
 $t^i$    & guaranteed               & possible (from \s{\trh(M_0)})  & possible (from \s{\rh(M_+)}) \\ \hline\hline
\s{t^{-5}}& \s{a^3q^3}               &                               &                       \\ \hline
\s{t^{-4}}& \s{aq^5}                 &                    &                   \\ \hline
\s{t^{-3}}& \s{a^3q^{-1}+aq^5 + a^2q^3(q^{N-4}+\dotsm +q^{-N+4}) }  & \s{aq^3} &    \\ \hline
\s{t^{-2}}& \s{aq+q^4(q^{N-3}+\dotsm+q^{-N+3})}  & \s{aq}    & \s{aq^3}            \\ \hline
\s{t^{-1}}& \s{aq+a^2q^{-1}(q^{N-4}+\dotsm + q^{-N+4})} & \s{a^{-1}q^3+aq^{-1}} &\s{aq} \\ \hline
\s{t^0}   & \s{a^{-1}q +aq^{-1}+ 3(q^{N-3}+\dotsm + q^{-N+3})} & \s{a^{-1}q+aq^{-3}} & \s{a^{-1}q^3+aq^{-1}} \\ \hline
\s{t}     & \s{a^{-2}q^2(q^{N-3}+\dotsm + q^{-N+3})}  & \s{a^{-1}q^{-1}} & \s{a^{-1}q+aq^{-3}}      \\ \hline
\s{t^2}   & \s{a^{-1}q^{-1} + q^{-4}(q^{N-3}+\dotsm + q^{-N+3})} & \s{a^{-1}q^{-3}}     & \s{a^{-1}q^{-1}}  \\ \hline
\s{t^3}   & \s{a^{-3}q + a^{-2}q^{-2}(q^{N-3}+\dotsm + q^{-N+3})}     &           & \s{a^{-1}q^{-3}}              \\ \hline
\s{t^4}   & \s{a^{-1}q^{-5}}                    &               &             \\ \hline
\s{t^5}   & \s{a^{-3}q^{-3}}                   &               &              \\ \hline
\end{tabular}\end{center}
\end{table}
\bigskip

The dimension of $E_\infty(1)$ has to be 1, living in homological degree 0. A straightforward computation shows that we already have this convergence in the column of guaranteed generators of $\rh(K_0)$. By inspection we see that if we promote the generator $t^{-3}aq^3$ in the column of possible generators from $\trh(M_0)$ than it would survive in $E_\infty(1)$. Therefore the exact sequence~\eqref{eq:ktotred} and the symmetry under $\psi$ imply that $t^{-3}aq^3$, $t^{-2}aq^3$, $t^2a^{-1}q^{-3}$ and $t^3a^{-1}q^{-3}$ must be discarded from the list of possible generators of $\rh(K_0)$. We present the updated list in Table~\ref{tab:2nd-KT}.
\begin{table}[h!]
\caption{Guaranteed and possible generators of $\rh(K_0)$ updated}
\label{tab:2nd-KT}
\newcommand\s[1]{\small{$#1$}}
\begin{center}
\begin{tabular}{|c|c|c|c|}
\hline
 $t^i$    & guaranteed               & possible (from \s{\trh(M_0)})  & possible (from \s{\rh(M_+)}) \\ \hline\hline
\s{t^{-5}}& \s{a^3q^3}               &                               &                       \\ \hline
\s{t^{-4}}& \s{aq^5}                 &                    &                   \\ \hline
\s{t^{-3}}& \s{a^3q^{-1}+aq^5 + a^2q^3(q^{N-4}+\dotsm +q^{-N+4}) }  &  &    \\ \hline
\s{t^{-2}}& \s{aq+q^4(q^{N-3}+\dotsm+q^{-N+3})}  & \s{aq}    &             \\ \hline
\s{t^{-1}}& \s{aq+a^2q^{-1}(q^{N-4}+\dotsm + q^{-N+4})} & \s{a^{-1}q^3+aq^{-1}} &\s{aq} \\ \hline
\s{t^0}   & \s{a^{-1}q +aq^{-1}+ 3(q^{N-3}+\dotsm + q^{-N+3})} & \s{a^{-1}q+aq^{-3}} & \s{a^{-1}q^3+aq^{-1}} \\ \hline
\s{t}     & \s{a^{-2}q^2(q^{N-3}+\dotsm + q^{-N+3})}  & \s{a^{-1}q^{-1}} & \s{a^{-1}q+aq^{-3}}      \\ \hline
\s{t^2}   & \s{a^{-1}q^{-1} + q^{-4}(q^{N-3}+\dotsm + q^{-N+3})} &           & \s{a^{-1}q^{-1}}  \\ \hline
\s{t^3}   & \s{a^{-3}q + a^{-2}q^{-2}(q^{N-3}+\dotsm + q^{-N+3})}     &             &              \\ \hline
\s{t^4}   & \s{a^{-1}q^{-5}}                    &               &             \\ \hline
\s{t^5}   & \s{a^{-3}q^{-3}}                   &               &              \\ \hline
\end{tabular}\end{center}
\end{table}
\bigskip

To determine whether the remaining possible generators are generators of $\rh(K_0)$ we use the spectral sequence $E_k(2)$.
The reduced Khovanov homology of $K_0$ (computed from KhoHo~\cite{khoho} with $q\ra q^{-1}$ to agree our conventions) is
\begin{gather}
\label{eq:KhK0}
\begin{split}
\overline{Kh}(K_0) &= 
q^{9}t^{-5} + q^{7}t^{-4} + q^{5}t^{-3}  + 2q^{3}t^{-2} + (q^{3} + q)t^{-1}  + 2(q+q^{-1})
\\ & \mspace{20mu}    
+ (q^{-1} + q^{-3})t + 2q^{-3}t^2 + q^{-5}t^3 + q^{-7}t^4 + q^{-9}t^5.
\end{split}
\end{gather}

To have $E_k(2)\Rightarrow\overline{Kh}(K_0)$ we have to promote the generators $t^{-2}aq$, $t^{-1}aq$, $t^{-1}a^{-1}q^3$, $a^{-1}q^3$, $aq^{-3}$, $taq^{-3}$, $ta^{-1}q^{-1}$ and $t^2a^{-1}q^{-1}$ (recall that in our conventions $d_k(N)=a^{-2k}q^{2Nk}t$).
A simple computation shows that $E_k(2)$ collapses after the second page \emph{i.e.} $E_2(2)=E_\infty(2)=\overline{Kh}(K_0)$ and that by promoting the four remaining generators the spectral sequence $E_k(2)$ would not converge to $\overline{Kh}(K_0)$. We leave the details to the reader.
After rearranging some terms we have that
\begin{gather}
\label{eq:HK0}
\begin{split}
\rh(K_0) &= 
   a^3q^3t^{-5} + aq^5t^{-4} + \bigl( a^3q^{-1}+a^2q^2[N-2] \bigr)t^{-3} + \bigl(2aq+q^4[N-2] \bigr)t^{-2} 
\\ & \mspace{20mu}
+ \bigl(a^2q^{-2}[N-2] + a^{-1}q^3+aq\bigr)t^{-1}
+ \bigl(a^{-1}q+aq^{-1}+a^{-1}q^3+aq^{-3} +3[N-2] \bigr)
\\ & \mspace{20mu}
+ \bigl(aq^{-3}+a^{-1}q^{-1}+ a^{-2}q^2[N-2]\bigr)t
+ \bigl(2a^{-1}q^{-1} + q^{-4}[N-2]\bigr)t^2
\\ & \mspace{20mu}
+ \bigl(a^{-3}q+a^{-2}q^{-2}[N-2]\bigr)t^3 + a^{-1}q^{-5}t^4 + a^{-3}q^{-3}t^5
.
\end{split}
\end{gather}
Using the exact sequence~\eqref{eq:totred} we obtain
\begin{gather}
\label{eq:HHK0}
\begin{split}
  \trh(K_0) &= \bigl( qt^{-1/2}+q^{-1}t^{1/2} \bigr) 
\bigl[a^3q^3t^{-5} + aq^5t^{-4} + a^3q^{-1}t^{-3} + 2aqt^{-2} + (a^{-1}q^3+aq)t^{-1} 
\\ & \mspace {-40mu}
+ (a^{-1}q^3+aq^{-3}) + (aq^{-3}+a^{-1}q^{-1})t + 2a^{-1}q^{-1}t^2 + a^{-3}qt^3 + a^{-1}q^{-5}t^4 + a^{-3}q^{-3}t^5 \bigr]
\\ & \mspace{-40mu}
+ qt^{-1/2}\bigl(a^3q^{-1}t^{-3} + aqt^{-2} + a^3q^{-5} + 2aq^{-3} + aq^{-1} + a^{-1}q^{-1} + aq^{-7} a^{-1}q^{-5}\bigr)
\\ & \mspace{-40mu}
+ q^{-1}t^{1/2}\bigl( aq^5 + a^{-1}q^7 + aq + 2a^{-1}q^3 + a^{-1}q + a^{-3}q^5 + a^{-1}q^{-1} + a^{-3}q \bigr)
.
\end{split}
\end{gather}
Finally, comparing $at^{-1/2}\trh(K_0)$ to $\rh(unknot)$ gives
\begin{gather*}
\label{eq:KT}
\begin{split}
\rh(K_\text{KT}) &=  
a^2q^4t^{-5} + (q^6 + a^2q^2)t^{-4} + (2a^2 + q^4)t^{-3} +  (3q^2 + a^2q^{-2} + q^4)t^{-2}
\\ & \mspace{-45mu}  
+ (a^{-2}q^6 +  a^2q^{-4} + a^{-2}q^4 + 2 + q^2)t^{-1} +
(3+3q^{-2} + a^{-2}q^2 + a^{-2}q^4)
\\ & \mspace{-45mu} 
+ (3a^{-2}q^2 + q^{-2} + 2a^{-2} + q^{-4})t 
+ (a^{-4}q^4 + q^{-6} + 2a^{-2} + q^{-4} + a^{-2}q^{-2})t^2 
\\ & \mspace{-45mu} 
+ (3a^{-2}q^{-2} + a^{-4}q^2 + a^{-2}q^{-4})t^3 + (2a^{-4} + a^{-2}q^{-4})t^4 
\\ & \mspace{-45mu}
+ (a^{-4}q^{-2} + a^{-2}q^{-6})t^5 + a^{-4}q^{-4}t^6.
\end{split}
\end{gather*}

\subsection{The Conway knot}\label{ssec:conway}

\begin{figure}[h!]
\centering
\figins{0}{1.8}{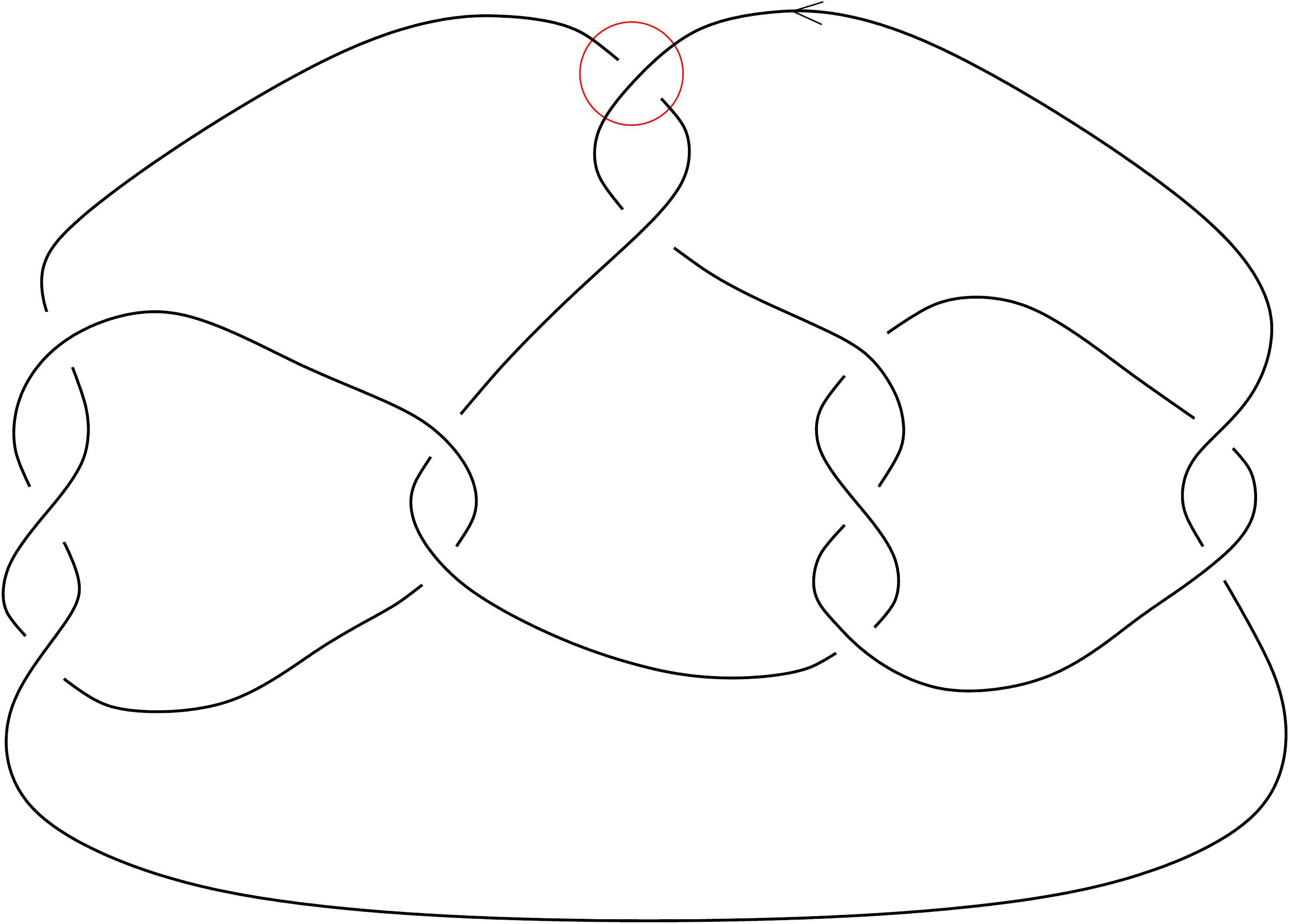}
\caption{The Conway knot}
\label{fig:conway}
\end{figure}

We follow the same method as in the calculation for the Kinoshita-Terasaka knot. A diagram of the Conway knot is given in Figure~\ref{fig:conway}.

\begin{figure}[h!]
\centering
\figins{0}{1.8}{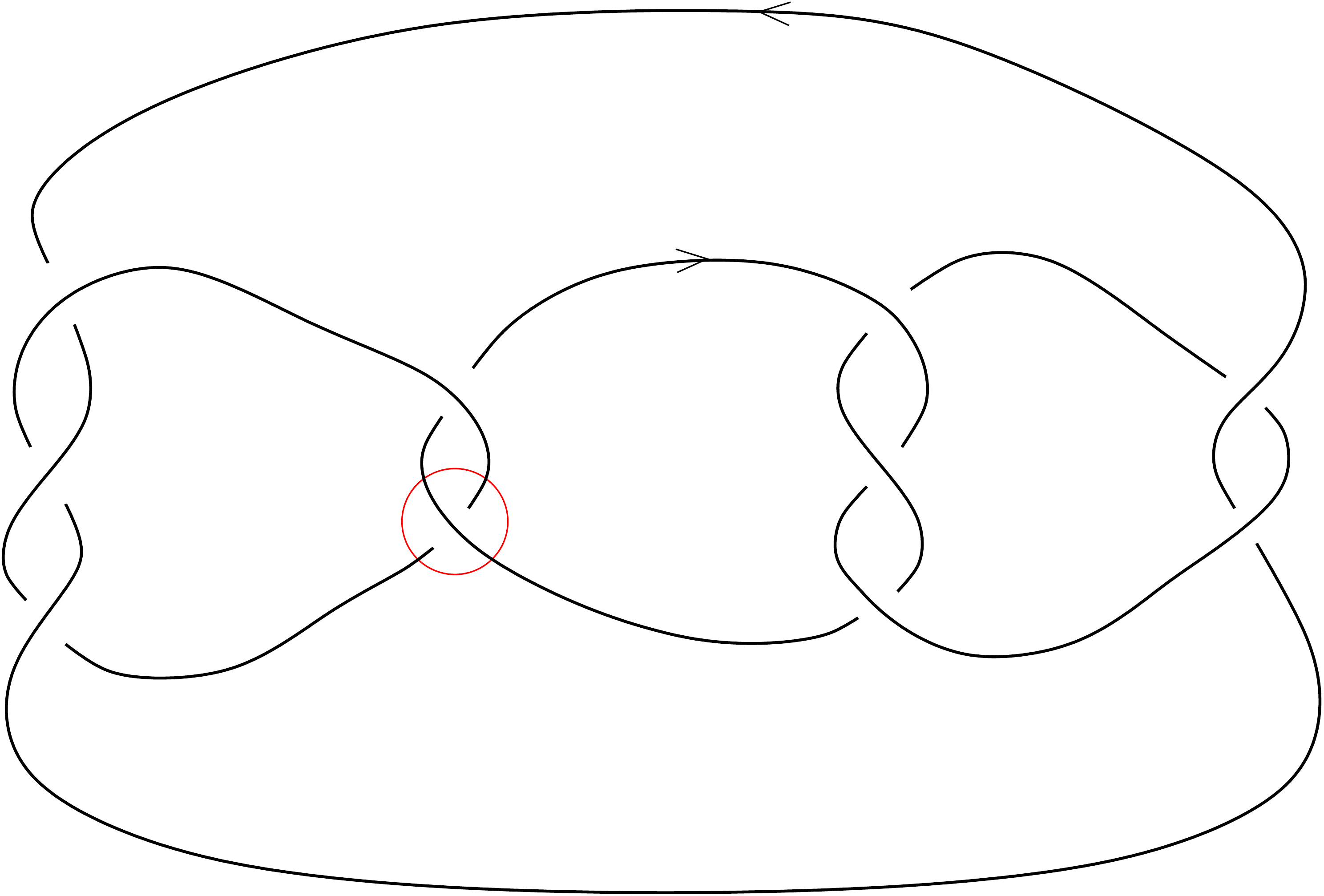}
\caption{Diagram $L_0$ obtained by taking the oriented resolution of the encircled crossing of the diagram of the Conway knot of Figure~\ref{fig:KT}. The (negative) encircled crossing is used for the long exact sequence~\eqref{eq:ktotred}}
\label{fig:cK0}
\end{figure}

Changing the (negative) encircled crossing we obtain the unknot, while resolving it results in the two-component link diagram $L_0$ of Figure~\ref{fig:cK0} and corresponds to the Pretzel link $P(3,-2,-3,2)$ which corresponds to the link $L_{10n59}$ in Thistlethwaite's table~\cite{katlas} and is amphicheiral. Using KhoHo~\cite{khoho} we find that
$$\overline{Kh}(L_{10n59})=\overline{Kh}(L_{10n36}),$$
with $\overline{Kh}(L_{10n36})$ given in Equation~\ref{eq:KhK0}.

\begin{figure}[h!]
\hair 2pt
\labellist
\pinlabel $N_+$ at  435 -50
\pinlabel $N_0$ at 1420 -50
\endlabellist
\centering
\figins{0}{1.75}{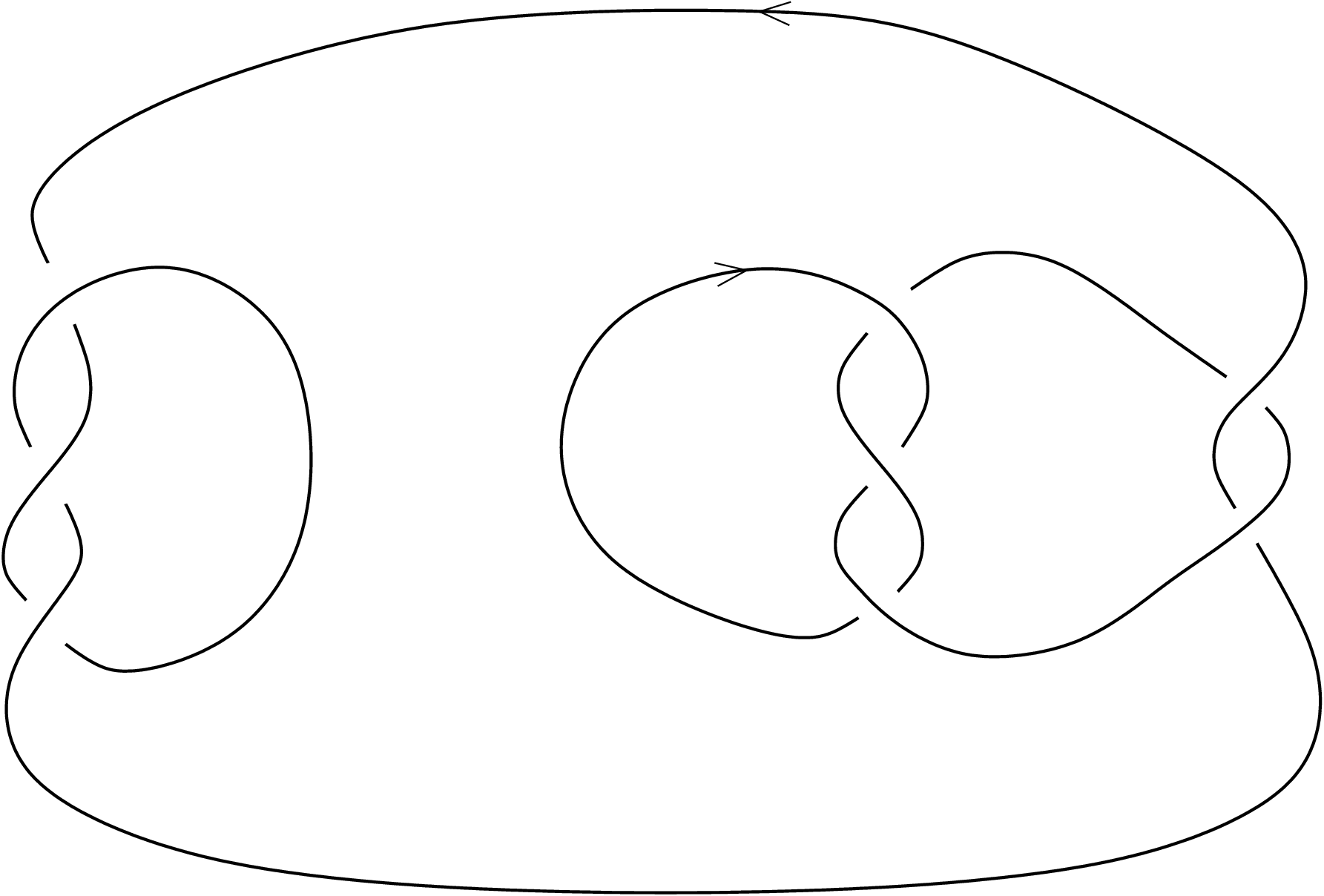}\qquad\ \figins{0}{1.75}{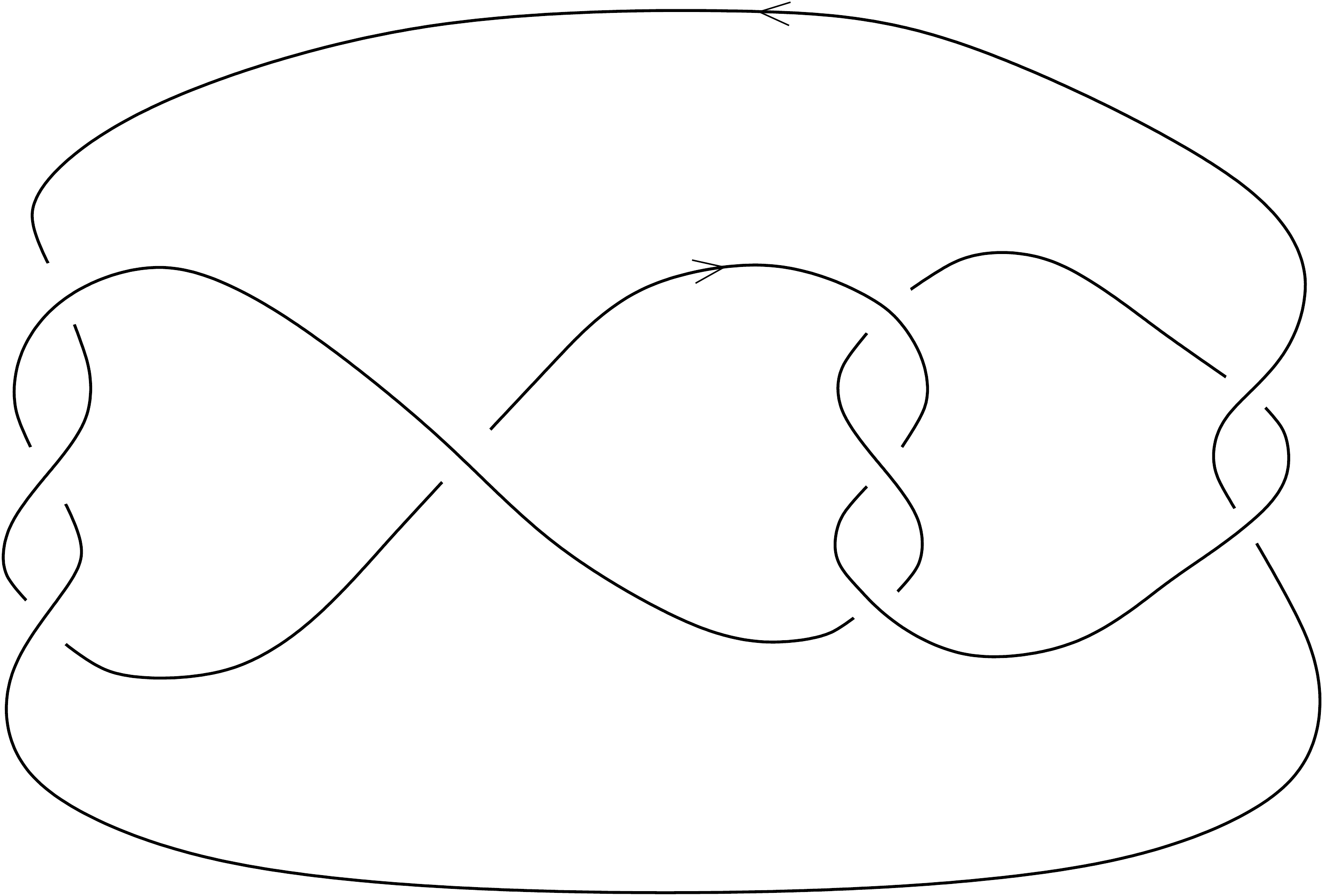}\\[2ex]
\caption{Diagrams $N_+$ and $N_0$ obtained by taking the oriented resolution of the encircled crossing of the diagram $L_0$ of Figure~\ref{fig:cK0} }
\label{fig:cM0}
\end{figure}

Changing and taking the oriented resolution of the encircled crossing of $L_0$ we obtain the diagrams $N_0$ and $N_+$ of Figure~\ref{fig:cM0}.
The diagram $N_+$ is isotopic to the connected sum of the (positive) Hopf link and the connected sum of positive trefoil and negative trefoil. The diagram $N_0$ is isotopic to the Pretzel knot $P(3,-3,2)$ which in turn corresponds to the mirror image of the knot $8_{20}$. The diagrams $N_0$ and $N_+$ are therefore isotopic to the diagrams $M_0$ and $M_+$ of Subsection~\ref{ssec:KT} respectively. This means that $\rh(L_0)=\rh(K_0)$ that is 
$$\rh\bigl(P(3,-2,-3,2)\bigr)=\rh\bigl(P(3,-2,2,-3)\bigr).$$

Since the other diagram obtained in the first step from the diagram for the Conway knot is the unknot (as in Subsection~\ref{ssec:KT}) we see that the HOMFLY-PT homologies of the Conway and Kinoshita-Terasaka knots are isomorphic and given by
\begin{gather*}
\label{eq:Hconway}
\begin{split}
\rh(K_\text{Conway}) =\rh(K_\text{KT}) &=  
a^2q^4t^{-5} + (q^6 + a^2q^2)t^{-4} + (2a^2 + q^4)t^{-3} +  (3q^2 + a^2q^{-2} + q^4)t^{-2}
\\ & \mspace{-45mu}  
+ (a^{-2}q^6 +  a^2q^{-4} + a^{-2}q^4 + 2 + q^2)t^{-1} +
(3+3q^{-2} + a^{-2}q^2 + a^{-2}q^4)
\\ & \mspace{-45mu} 
+ (3a^{-2}q^2 + q^{-2} + 2a^{-2} + q^{-4})t 
+ (a^{-4}q^4 + q^{-6} + 2a^{-2} + q^{-4} + a^{-2}q^{-2})t^2 
\\ & \mspace{-45mu} 
+ (3a^{-2}q^{-2} + a^{-4}q^2 + a^{-2}q^{-4})t^3 + (2a^{-4} + a^{-2}q^{-4})t^4 
\\ & \mspace{-45mu}
+ (a^{-4}q^{-2} + a^{-2}q^{-6})t^5 + a^{-4}q^{-4}t^6.
\end{split}
\end{gather*}

Note that $d_k(N)=0$ for $N\geq 3$ for $\rh(K_{\text{Conway}})$ and $\rh(K_{\text{KT}})$ for degree reasons. This implies that $\rhn(K_{\text{Conway}})=\rhn(K_{\text{KT}})$ for all $N\geq 3$. We already knew that the same holds for $N=2$ by direct computation. Since the knot Floer homology of these two knots differ~\cite{OS}, this fact might be interesting for someone trying to find a relation (e.g. a spectral sequence) between Khovanov-Rozansky homology and knot Floer homology.


\vspace*{1cm}

\noindent {\bf Acknowledgements} 
We thank Jacob Rasmussen for the enlightening exchanges of email about the topic of this paper.

The authors were supported by the 
Funda\c {c}\~{a}o para a Ci\^{e}ncia e a Tecnologia (ISR/IST plurianual funding) through the
programme ``Programa Operacional Ci\^{e}ncia, Tecnologia, Inova\-\c
{c}\~{a}o'' (POCTI) and the POS Conhecimento programme, cofinanced by the European Community 
fund FEDER.


\end{document}